%&amstex
\input amstex
\input amsppt.sty
\magnification=\magstep1
\hsize=30truecc
\baselineskip=16truept
\vsize=22.2truecm
%\NoBlackBoxes
\nologo
\pageno=1
\topmatter
\TagsOnRight

\def\Z{\Bbb Z}

\def\Q{\Bbb Q}

\def\al{\alpha}
\def\l{\left}
\def\r{\right}
\def\bg{\bigg}
\def\({\bg(}
\def\){\bg)}
\def\[{\bg[}
\def\]{\bg]}
\def\t{\text}
\def\f{\frac}
\def\em{\emptyset}
\def\se {\subseteq}

\def\bi{\binom}

\def\ls{\leqslant}
\def\gs{\geqslant}
\def\mo{\t{\rm mod}}

\def\Proof{\noindent{\it Proof}}

\def\Remark{\medskip\noindent{\it  Remark}}
\def\Ack{\medskip\noindent {\bf Acknowledgment}}
\hbox{Final version for Proc. Amer. Math. Soc.}
\bigskip
\title A sharp result on $m$-covers\endtitle
\author Hao Pan and Zhi-Wei Sun\endauthor
\address Department of Mathematics, Nanjing University,
Nanjing 210093, People's Republic of China
\newline\indent
{\it E-mail addresses}: {\tt haopan79\@yahoo.com.cn},\ \ \tt{zwsun\@nju.edu.cn}
\endaddress

\abstract  Let $A=\{a_s+n_s\Z\}_{s=1}^k$ be a finite system of arithmetic
sequences which forms an $m$-cover of $\Z$ (i.e., every integer
belongs at least to $m$ members of $A$). In this paper we show the
following sharp result: For any positive integers $m_1,\ldots,m_k$
and $\theta\in[0,1)$, if there is $I\se\{1,\ldots,k\}$ such
that the fractional part of $\sum_{s\in I} m_s/n_s$ is $\theta$,
then there are at least $2^m$ such subsets of $\{1,\ldots,k\}$.
This extends an earlier result of M. Z. Zhang and an extension by Z. W.
Sun. Also, we generalize the above result to $m$-covers of the
integral ring of any algebraic number field with a power integral basis.
\endabstract

\thanks 2000 {\it Mathematics Subject Classifications}:\,Primary 11B25;
Secondary 11B75, 11D68, 11R04.
\newline\indent The second author is responsible for
communications, and supported by the National Science Fund for
Distinguished Young Scholars (No. 10425103) in China.
\endthanks
\endtopmatter
\document

\heading{1. Introduction}\endheading

For an integer $a$ and a positive integer $n$, we simply let
$a(n)$ represent the set $a+n\Z=\{x\in\Z:\, x\equiv a\ (\mo\ n)\}$.
Following Sun [S95, S96] we call a finite system
$$A=\{a_s(n_s)\}_{s=1}^k\tag 1.1$$
of such sets an {\it $m$-cover} of $\Z$ (where
$m\in\{1,2,3,\ldots\})$ if every integer lies in at least $m$
members of $(1.1)$. We use the term {\it cover} (or covering system)
instead of $1$-cover. For problems and results in this area, the
reader may consult [G04, pp.\,383--390], [PS] and [S05].
P. Erd\H os [E97] once said:
 {\it ``Perhaps my favorite problem of all concerns covering systems}."
\smallskip

{\it Example} 1.1. For each integer $m\gs1$,
there is an $m$-cover of $\Z$ which is not the union of two covers of $\Z$.
To wit, we let
$p_1,\ldots,p_r$ be distinct primes with
$r\gs 2m-1$, and set $N=p_1\cdots p_r$.
Clearly $A_*=\{\prod_{s\in I}p_s(N)\}_{I\se\{1,\ldots,r\},\,|I|\gs m}$
does not cover any integer relatively prime to $N$.
Let $a_1,\ldots,a_n$ be the list
of those integers in $\{0,1,\ldots,N-1\}$
not covered by $A_*$ with each occurring exactly $m$ times.
If $x\in\Z$ is covered by $A_*$, then $x\in\bigcap_{s\in I}0(p_s)$
for some $I\se\{1,\ldots,r\}$ with $|I|\gs m$. Therefore
$$A_0=\{0(p_1),\ldots,0(p_r),a_1(N),\ldots,a_n(N)\}$$
forms an $m$-cover of $\Z$. Suppose that
$I_1\cup I_2=\{1,\ldots,r\}$,
$J_1\cup J_2=\{1,\ldots,n\}$ and $I_1\cap I_2=J_1\cap J_2=\em$.
For $i=1,2$ let $A_i$ be the system consisting of
those $0(p_s)$ with $s\in I_i$ and those $a_t(N)$ with $t\in J_i$.
We claim that  $A_1$ or $A_2$ is not a cover of $\Z$.
Without loss of generality, let us assume that $|I_1|\ls |I_2|$.
Since $2|I_2|\gs|I_1|+|I_2|>2(m-1)$, we have $|I_2|\gs m$ and hence
$\prod_{s\in I_2}p_s$ is covered by $A_*$.
Therefore $\prod_{s\in I_2}p_s\not\in\bigcup_{t=1}^na_t(N)$.
Clearly
$\prod_{s\in I_2}p_s$ is not covered by $\{0(p_s)\}_{s\in I_1}$ either.
Thus $A_1$ does not form a cover of $\Z$.
\medskip

 By means of the Riemann zeta function, in 1989 M. Z. Zhang [Z89]
proved that if $(1.1)$ forms a cover of $\Z$ then
$\sum_{s\in I}1/n_s$ is a positive integer for some $I\subseteq\{1,\ldots,k\}$.

 Let $m_1,\ldots,m_k$ be any positive integers. If $(1.1)$ is a cover
of $\Z$, then $\{a_s+(n_s/m_s)\Z\}_{s=1}^k$ is also a cover of
$\Z$ and hence Theorem 2 of [S95] indicates that for any
$J\subseteq\{1,\ldots,k\}$ there is an $I\subseteq\{1,\ldots,k\}$
with $I\not=J$ such that $\{\sum_{s\in I}m_s/n_s\}=\{\sum_{s\in
J}m_s/n_s\}$, where $\{\alpha\}$ denotes the fractional part of
a real number $\alpha$. When $J=\em$ and $m_1=\cdots=m_k=1$, this yields
Zhang's result. In 1999 Z. W. Sun [S99] proved further that
if $(1.1)$ forms an $m$-cover of $\Z$ then for any
$J\subseteq\{1,\ldots,k\}$ we have
$$\bg|\bg\{I\subseteq\{1,\ldots,k\}:\, I\not=J\ \text{and}\
\bg\{\sum_{s\in I}\frac{m_s}{n_s}\bg\}
=\bg\{\sum_{s\in J}\frac{m_s}{n_s}\bg\}\bg\}\bg|\gs m.$$

 In this paper we will show the following sharp result.
\proclaim{Theorem 1.1} Let $(1.1)$ be an $m$-cover of $\Z$, and let
$m_1,\ldots,m_k$ be any integers. Then for any $0\ls\theta<1$ the
set
$$I_A(\theta)=\bg\{I\subseteq\{1,\ldots,k\}:\,
\bg\{\sum_{s\in I}\frac{m_s}{n_s}\bg\}=\theta\bg\}\tag1.2$$
has at least $2^m$ elements if it is nonempty.
\endproclaim
\Remark\ 1.1. Clearly $m$ copies of $0(1)$ form an $m$-cover of
$\Z$. This shows that the lower bound in Theorem 1.1 is best
possible. \proclaim{Corollary 1.1} Let $(1.1)$ be an $m$-cover of
$\Z$, and let $m_1,\ldots,m_k$ be any integers. Then
$|S(A)|\ls 2^{k-m}$ where
$$S(A)=\bg\{\bg\{\sum_{s\in I}\frac{m_s}{n_s}\bg\}:\,
I\subseteq\{1,\ldots,k\}\bg\}.\tag1.3$$
\endproclaim
\Proof. As $|I_A(\theta)|\gs 2^m$ for all $\theta\in S(A)$, we
have
$$|S(A)|2^m\ls |\{I:\,I\subseteq\{1,\ldots,k\}\}|=2^k$$
and hence $|S(A)|\ls 2^{k-m}$.\qed

\Remark\ 1.2. Sun [S95, S96] showed that if $m_1,\ldots,m_k$ are
relatively prime to $n_1,\ldots,n_k$ respectively then $(1.1)$
forms an $m$-cover of $\Z$ whenever it covers $|S(A)|$ consecutive
integers at least $m$ times.
\medskip

\proclaim{Corollary 1.2} Suppose that $(1.1)$
forms an $m$-cover of $\Z$ but $\{a_s(n_s)\}_{s=1}^{k-1}$ does not.
If the covering function
$w_A(x)=|\{1\ls s\ls k:\,x\in a_s(n_s)\}|$ is periodic modulo $n_k$, then
for any $r=0,\ldots,n_k-1$ we have
$$\bg|\bg\{I\se\{1,\ldots,k-1\}:\,
\bg\{\sum_{s\in I}\f1{n_s}\bg\}=\f r{n_k}\bg\}\bg|\gs 2^{m-1}.\tag1.4$$
\endproclaim
\Proof. By Theorem 1 of Sun [S06],
$$\bg|\bg\{\bg\lfloor\sum_{s\in I}\f1{n_s}\bg\rfloor:\
I\se\{1,\ldots,k-1\}\ \t{and}\ \bg\{\sum_{s\in I}\f1{n_s}\bg\}=\f r{n_k}\bg\}\bg|\gs m.$$
In particular,
$\{\sum_{s\in I}1/n_s\}=r/n_k$ for some
 $I\se\{1,\ldots,k-1\}$, and hence (1.4) holds in the case $m=1$.
For $A_k=\{a_s(n_s)\}_{s=1}^{k-1}$, clearly $w_{A_k}(x)\gs m-1$ for all $x\in\Z$.
In the case $m>1$, we obtain (1.4)
by applying Theorem 1.1 to $A_k$ with  $m_1=\cdots=m_{k-1}=1$
and $\theta=r/n_k$. \qed

\Remark\ 1.3. When $n_k$ is divisible by all the moduli $n_1,\ldots,n_k$, Corollary 1.2
was stated by the second author in [S03, Theorem 2.5].
When $w_A(x)=m$ for all $x\in\Z$, the following result stronger than (1.4) was proved in [S97]:
$$\bg|\bg\{I\se\{1,\ldots,k-1\}:\,\sum_{s\in I}\f1{n_s}=n+\f r{n_k}\bg\}\bg|\gs\bi{m-1}n$$
for every $n=0,\ldots,m-1$.
\medskip

For an algebraic number field $K$, let $O_K$ be the ring of
algebraic integers in $K$. For $\alpha, \beta\in O_K$,
we set
$$\alpha+\beta O_K=\{\al+\beta\omega:\, \omega\in O_K\}$$
and call it a residue class in $O_K$.
For a finite system
$$\Cal A=\{\alpha_s+\beta_s O_K\}_{s=1}^k\tag1.5$$
of such residue classes,
if $|\{1\ls s\ls k:\,x\in\al_s+\beta_s O_K\}|\gs m$
for all $x\in O_K$ (where $m\in\{1,2,3,\ldots\}$), then we call $\Cal A$
an $m$-cover of $O_K$.
Covers of the ring $\Z[\sqrt{-2}]=O_{\Q(\sqrt{-2})}$ were
investigated by J. H. Jordan [J68].

An algebraic number field $K$ of degree $n$ is said to have a {\it power integral basis}
if there is $\gamma\in O_K$ such that $1,\gamma,\ldots,\gamma^{n-1}$ form a basis of $O_K$ over $\Z$.
It is well known that all quadratic fields and
cyclotomic fields have power integral bases.

 Here is a generalization of Theorem 1.1.

\proclaim{Theorem 1.2} Let $K$ be an algebraic number field with
a power integral basis. Suppose that $(1.5)$ forms
an $m$-cover of $O_K$, and let $\omega_1,\ldots,\omega_k\in O_K$. Then, for any $\mu\in K$, the set
$$\bg\{I\subseteq\{1,\ldots,k\}:\,
\sum_{s\in I}\f{\omega_s}{\beta_{s}}\in\mu+O_K\bg\}$$ is empty or it has at
least $2^m$ elements.
\endproclaim

\Remark\ 1.4. We conjecture that the requirement in Theorem 1.2 that $K$ has a power integral basis
can be cancelled.

\heading 2. Proof of Theorem 1.1\endheading

\proclaim{Lemma 2.1} Let $(1.1)$ be an $m$-cover of $\Z$. Let
$m_1,\ldots,m_k\in\Z$ and define $S(A)$ as in $(1.3)$. Then, for any given
$\theta\in S(A)$, there exists $t\in\{1,\ldots,k\}$ such that
both $\theta$ and $\{\theta-m_t/n_t\}$ lie in $S(A_t)$, where
$A_t=\{a_s(n_s)\}_{1\ls s\ls k,\ s\not=t}$.
\endproclaim
\Proof. Choose a maximal $J\subseteq\{1,\ldots,k\}$ such that
$\{\sum_{s\in J}m_s/n_s\}=\theta$. As $(1.1)$ is a cover of $\Z$, by
[S95, Theorem 2] or [S99, Theorem 1(i)] there exists an
$I\subseteq\{1,\ldots,k\}$ for which $I\not=J$ and
$\{\sum_{s\in I}m_s/n_s\}=\theta$. Note that $J\not\subseteq I$
and hence $t\in J\setminus I$ for some $1\ls t\ls k$. Clearly
$\theta=\{\sum_{s\in I}m_s/n_s\}\in S(A_t)$ and also
$\{\theta-m_t/n_t\}=\{\sum_{s\in J\setminus\{t\}}m_s/n_s\}\in
S(A_t)$. This concludes the proof.\qed

\medskip\noindent {\it Proof of Theorem 1.1}.
We use induction on $m$.

The $m=1$ case, as mentioned above, has been handled in [S95, S99].

Now let $m>1$ and assume that Theorem 1.1 holds for smaller
positive integers. Let $\theta\in S(A)$. In light of Lemma 2.1,
there is a $t\in\{1,\ldots,k\}$ such that both
$\theta$ and $\theta'=\{\theta-m_t/n_t\}$ lie in $S(A_t)$. As
$A_t$ forms a $(m-1)$-cover of $\Z$,  by the induction hypothesis
we have $|I_{A_t}(\theta)|\gs 2^{m-1}$
and $|I_{A_t}(\theta')|\gs 2^{m-1}$. Observe that
$$I_A(\theta)=I_{A_t}(\theta)\cup\{I\cup\{t\}:
\, I\in I_{A_t}(\theta')\}.$$
Therefore
$$|I_A(\theta)|=|I_{A_t}(\theta)|+|I_{A_t}(\theta')|
\gs2^{m-1}+2^{m-1}=2^m.$$
We are done.\qed

\heading {3. Proof of Theorem 1.2}\endheading

At first we give a lemma on algebraic number fields with power integral bases.

\proclaim{Lemma 3.1} Let $K$ be an algebraic number field with a power integral
basis $1,\gamma,\ldots,\gamma^{n-1}$.
For any $\mu=\sum_{r=0}^{n-1}\mu_r\gamma^r\in K$ with $\mu_0,\ldots,\mu_{n-1}\in\Q$,
we have
$$\mu\in O_K\iff \psi(\mu),\psi(\mu\gamma),\ldots,\psi(\mu\gamma^{n-1})\in\Z,$$
where $\psi(\mu)$ denotes the last coordinate $\mu_{n-1}$.
\endproclaim
\Proof. If $\mu\in O_K$, then $\mu,\mu\gamma,\ldots,\mu\gamma^{n-1}\in O_K$
and hence $\psi(\mu\gamma^j)\in\Z$ for every $j=0,\ldots,n-1$.

Now assume that $\psi(\mu\gamma^j)\in\Z$ for all $j=0,\ldots,n-1$.
We want to show that $\mu\in O_K$ (i.e., $\mu_0,\ldots,\mu_{n-1}\in\Z$).
Clearly $\mu_{n-1}=\psi(\mu\gamma^0)\in\Z$. If
$0\ls r<n-1$ and $\mu_{r+1},\ldots,\mu_{n-1}\in\Z$, then
$$\align\mu_r=&\psi(\mu_0\gamma^{n-1-r}+\mu_1\gamma^{n-r}+\cdots+\mu_r\gamma^{n-1})
\\=&\psi(\mu\gamma^{n-1-r})-\psi(\mu_{r+1}\gamma^n+\cdots+\mu_{n-1}\gamma^{2n-2-r})
\endalign$$
and hence $\mu_r\in\Z$ since $\mu_{r+1}\gamma^n+\cdots+\mu_{n-1}\gamma^{2n-2-r}\in O_K$.
So, by induction, $\mu_r\in\Z$ for all $r=0,\ldots,n-1$. We are done. \qed

\medskip
\noindent{\it Proof of Theorem 1.2}. In the spirit of the proof
of Theorem 1.1, it suffices to handle the case $m=1$.
That is, we should prove that for any $J\se\{1,\ldots,k\}$ there is $I\se\{1,\ldots,k\}$
with $I\not=J$ such that $\sum_{s\in I}\omega_s/\beta_s-\sum_{s\in J}\omega_s/\beta_s\in O_K$.

Let $\{1,\gamma,\ldots,\gamma^{n-1}\}$ be a power integral basis of $K$,
and define $\psi$ as in Lemma 3.1.

Let $x_0,\ldots,x_{n-1}\in\Z$ and $x=\sum_{r=0}^{n-1}x_r\gamma^r$.
Since $-x\in O_K$
is covered by $\Cal A=\{\al_s+\beta_s O_K\}_{s=1}^k$, we have
$$\align 0=&\prod_{s=1}^k\l(1-e^{2\pi i\psi(\omega_s(x+\alpha_s)/\beta_s)}\r)
=\sum_{I\subseteq\{1,\ldots,k\}}(-1)^{|I|}\prod_{s\in I}
e^{2\pi i\psi(\omega_s(x+\alpha_s)/\beta_s)}
\\=&\sum_{I\subseteq\{1,\ldots,k\}}(-1)^{|I|}\prod_{s\in I}
e^{2\pi i(\psi(\omega_s\alpha_s/\beta_s)+\sum_{r=0}^{n-1}x_r\psi(\omega_s\gamma^{r}/\beta_s))}
\\=&\sum_{I\subseteq\{1,\ldots,k\}}(-1)^{|I|}e^{2\pi i\psi(\sum_{s\in
I}\omega_s\alpha_s/\beta_s)}\prod_{r=0}^{n-1}
e^{2\pi ix_r\psi(\sum_{s\in I}\omega_s\gamma^{r}/\beta_s)}
\\=&\sum_{\theta_0\in S_0}e^{2\pi ix_0\theta_0}
\sum_{\theta_1\in S_1}e^{2\pi ix_1\theta_1}\cdots\sum_{\theta_{n-1}\in S_{n-1}}
e^{2\pi ix_{n-1}\theta_{n-1}}f(\theta_0,\ldots,\theta_{n-1}),
\endalign$$
where
$$S_r=\bg\{\bg\{\psi\(\sum_{s\in I}\f{\omega_s\gamma^r}{\beta_s}\)\bg\}:\, I\se\{1,\ldots,k\}\bg\}$$
and
$$f(\theta_0,\ldots,\theta_{n-1})=\sum\Sb I\se\{1,\ldots,k\}
\\\{\psi(\sum_{s\in I}\omega_s\gamma^r/\beta_s)\}=\theta_r
\\ \t{for all}\ r=0,\ldots,n-1\endSb
(-1)^{|I|}e^{2\pi i\psi(\sum_{s\in I}\omega_s\al_s/\beta_s)}.$$

 For each $r=0,\ldots,n-1$, if $\sum_{\theta_r\in S_r}e^{2\pi ix_r\theta_r}F(\theta_r)=0$
for all $x_r=0,\ldots,|S_r|-1$, then $F(\theta_r)=0$ for every $\theta_r\in S_r$, because
the Vandermonde determinant $\det(e^{2\pi i x_r\theta_r})_{0\ls x_r<|S_r|,\,\theta_r\in S_r}$
does not vanish. So, by the above, we have
$f(\theta_0,\ldots,\theta_{n-1})=0$ for all $\theta_0\in S_0,\ldots,\theta_{n-1}\in S_{n-1}$.

 Now suppose that $\mu\in K$ and $\sum_{s\in
 J}\omega_s/\beta_s\in\mu+O_K$ for a unique subset $J$ of
 $\{1,\ldots,k\}$. We want to deduce a contradiction.

 Set $\theta_r=\{\psi(\mu\gamma^r)\}$ for $r=0,\ldots,n-1$.
 For any $I\se\{1,\ldots,k\}$ we have
 $$\align&\bg\{\psi\(\sum_{s\in I}\f{\omega_s\gamma^r}{\beta_s}\)\bg\}=\theta_r\ \ \t{for all}\
 r=0,\ldots,n-1
 \\\iff&\psi\(\(\sum_{s\in I}\f{\omega_s}{\beta_s}-\mu\)\gamma^r\)\in\Z\ \ \t{for all}\
 r=0,\ldots,n-1
\\\iff&\sum_{s\in I}\f{\omega_s}{\beta_s}\in\mu+O_K\qquad\t{(by Lemma 3.1)}
\\\iff&I=J.
\endalign$$
Thus the expression of $f(\theta_0,\ldots,\theta_{n-1})$ only contains one summand, and therefore
$$0=f(\theta_0,\ldots,\theta_{n-1})=(-1)^{|J|}e^{2\pi i\psi(\sum_{s\in J}\omega_s\al_s/\beta_s)}\not=0$$
which is a contradiction.

 The proof of Theorem 1.2 is now complete. \qed

 \Ack. The authors are indebted to the referee for his/her valuable suggestions.

\enddocument